\documentclass[11pt,a4paper,twoside]{article}
\usepackage{amssymb}
\usepackage{amsmath}
\usepackage{amsthm}
\usepackage{rotate}
\usepackage{graphicx}

\usepackage[T1]{fontenc}
   \usepackage{fancyhdr}
\usepackage{afterpage}

\usepackage{latexcad}
   \renewcommand{\headrulewidth}{0.4pt}
\newtheorem{theorem}{Theorem}

\newtheorem{lemma}[theorem]{Lemma}

\newtheorem{proposition}[theorem]{Proposition}

\theoremstyle{definition}

\begin{document}

\def\vm#1#2#3#4{\left(\begin{array}{cc}#1 & #2\\ #3 & #4 \end{array}\right)}
\def\paige#1{M^*(#1)}
\def\cyclic#1{C_{#1}}
\def\dpr#1#2{#1\cdot#2}
\def\vpr#1#2{#1\times#2}
\def\span#1{\langle#1\rangle}
\def\field#1{GF(#1)}
\def\TripleCoxGroup#1#2#3{(#1,\,#2,\, #3)}
\def\QuadrupleCoxGroup#1#2#3#4{(#1,\,#2\ |\ #3,\,#4)}
\def\CoxGroup#1{\QuadrupleCoxGroup{3}{3}{3}{#1}}
\def\lattice#1{\mathcal{L}(#1)}
\def\image#1{\mathrm{Im} #1}
\def\neutral{E}
\def\gcd#1#2{\mathrm{gcd}(#1,\,#2)}
\def\cdeg#1#2{\mathrm{deg}_{#1}(#2)}
\def\order#1{o(#1)}
\def\aut#1{\mathrm{Aut}(#1)}
\def\join{\vee}

\afterpage{\fancyhead[CO]{$\rule{31mm}{0mm}$\small The abstract
groups $(3, 3\, |\, 3, p)$ \hfill\thepage}
\fancyhead[RE]{\thepage\hfill\small P. Vojt\v echovsk\'y
$\rule{48mm}{0mm}$} }

\begin{center}
\vspace*{2pt}
{\Large \textbf{The abstract groups (3, 3 | 3, p),}}\\[3mm]
{\Large\textbf{their subgroup structure,}}\\[3mm]
{\Large\textbf{and their significance for Paige loops}}\\[36pt]
{\large \textsf{\emph{Petr Vojt\v echovsk\'y}}}
\\[36pt]
\textbf{Abstract}
\end{center}
{\footnotesize{For most (and possibly all) non-associative finite
simple Moufang loops, three generators of order $3$ can be chosen
so that each two of them generate a group isomorphic to $(3, 3\,
|\, 3, p)$. The subgroup structure of $(3, 3\, |\, 3, p)$ depends
on the solvability of a certain quadratic congruence, and it is
described here in terms of generators.}}

 \footnote{group $\CoxGroup{p}$, loop generator, quadratic congruence}

\section*{\centerline{1. Introduction}}
\rm
 Moufang loops and, more generally, diassociative loops are
usually an abundant source of two-generated groups. In the end,
this is what diassociativity is all about: every two elements
generate an associative subloop, i.e.  a group. (We refer the
reader not familiar with the theory of loops to
\cite{Pflugfelder}.) This short paper emerged as an offshoot of
our larger-scale program to fully describe the subloop structure
of all non-associative finite simple Moufang loops, sometimes
called \emph{Paige loops}.

Let $\paige{q}$ denote the Paige loop constructed over $F=GF(q)$ as in
\cite{Paige}. That is, $\paige{q}$ consists of vector matrices
\begin{displaymath}
    M=\vm{a}{\alpha}{\beta}{b},
\end{displaymath}
where $a$, $b\in F$, $\alpha$, $\beta\in F^{3}$,
$\det{M}=ab-\dpr{\alpha}{\beta}=1$, and where $M$ is identified with $-M$. The
multiplication in $\paige{q}$ coincides with the Zorn matrix multiplication
\begin{displaymath}
    \vm{a}{\alpha}{\beta}{b}\vm{c}{\gamma}{\delta}{d}
    =\vm{ac+\dpr{\alpha}{\delta}}
        {a\gamma+d\alpha-\vpr{\beta}{\delta}}
        {c\beta+b\delta+\vpr{\alpha}{\gamma}}
        {\dpr{\beta}{\gamma}+bd},
\end{displaymath}
where $\dpr{\alpha}{\beta}$ (resp. $\vpr{\alpha}{\beta}$) is the standard dot
product (resp. cross product).

We have shown in \cite[Theorem 1.1]{PevoLondon} that every $\paige{q}$ is
three-generated, and when $q\ne 9$ is odd or $q=2$ then the generators can be
chosen as
\begin{equation}\label{Eq:OddPaigeGensAgain}
    g_1=\vm{1}{e_1}{0}{1},\;\;g_2=\vm{1}{e_2}{0}{1},\;\;
        g_3=\vm{0}{ue_3}{-u^{-1}e_3}{1},
\end{equation}
where $u$ is a primitive element of $F$ (cf. \cite[Proposition
4.1]{PevoLondon}). In particular, note that $g_1$, $g_2$ and $g_3$ generate
$\paige{p}$ for every prime $p$. We find it more convenient to use another set
of generators.

\begin{proposition} Let $q\ne 9$ be an odd
prime power or $q=2$. Then $\paige{q}$ is generated by three
elements of order three.

\begin{proof} Let us introduce
\begin{eqnarray*}
    g_4=g_3g_1&=&\vm{0}{(0,0,u)}{(0,u,-u^{-1})}{1},\\
    g_5=g_3g_2&=&\vm{0}{(0,0,u)}{(-u,0,-u^{-1})}{1}.
\end{eqnarray*}
It follows from $(\ref{Eq:OddPaigeGensAgain})$ that $\paige{q}$ is
generated by $g_3$, $g_4$, and $g_5$. One easily verifies that
these elements are of order $3$.
\end{proof}
\end{proposition}

The groups $\span{g_3,\,g_4}$, $\span{g_3,\,g_5}$ and
$\span{g_4,\,g_5}$ play therefore a prominent role in the lattice
of subloops of $\paige{q}$. As we prove in Section 3, each of them
is isomorphic to the group $\CoxGroup{p}$, defined below.

\section*{\centerline{2. The abstract groups (3, 3 | 3, p)}}

 The two-generated abstract groups $\QuadrupleCoxGroup{l}{m}{n}{k}$
defined by presentations
\begin{equation}
    \QuadrupleCoxGroup{l}{m}{n}{k}
        =\span{x,\,y\ |\ x^l=y^m=(xy)^n=(x^{-1}y)^k}\label{Eq:CoxQuadruple}
\end{equation}
were first studied by Edington \cite{Edington}, for some small values of $l$,
$m$, $n$ and $k$. The notation we use was devised by Coxeter \cite{Coxeter39}
and Moser \cite{CoxeterMoser}, and has a deeper meaning that we will not
discuss here. From now on, we will always refer to presentation
$(\ref{Eq:CoxQuadruple})$ when speaking about $\QuadrupleCoxGroup{l}{m}{n}{k}$.

The starting point for our discussion is Theorem \ref{Th:Edington}, due to
Edington \cite[Theorem IV and pp. 208--210]{Edington}. (Notice that there is a
typo concerning the order of $\CoxGroup{n}$, and a misprint claiming that
$\CoxGroup{3}$ is isomorphic to $A_4$.). For the convenience of the reader, we
give a short, contemporary proof.

\begin{theorem}[Edington]\label{Th:Edington}
The group $G=\CoxGroup{n}$ exists for every $n\geqslant 1$, is of
order $3n^2$, and is non-abelian when $n>1$. It contains a normal
subgroup $H=\span{x^2y,\,xy^2}\cong\cyclic{n}\times\cyclic{n}$. In
particular, $G\cong\cyclic{3}$ when $n=1$, $G\cong A_4$ when
$n=2$, and $G$ is the unique non-abelian group of order $27$ and
exponent $3$ when $n=3$.
\begin{proof}
Verify that $\CoxGroup{1}$ is isomorphic to $\cyclic{3}$. Let
$n>1$. Since $x(x^2y)x^{-1}=yx^{-1}=y(x^2y)y^{-1}\in H$, and
$x^{-1}(xy^2)x=y^2x=y^{-1}(xy^2)y\in H$, the subgroup $H$ is
normal in $G$. It is an abelian group of order at most $n^2$ since
$x^2y\cdot xy^2=x(xy)^2y=x(xy)^{-1}y=xy^2\cdot x^2y$. Clearly,
$G/H\cong\cyclic{3}$ (enumeration of cosets works fine), and hence
$|G|=3|H|\leqslant 3n^2$.

Let $N=\span{a}\times\span{b}\cong\cyclic{n}\times\cyclic{n}$, and
$K=\span{f}\leqslant\aut{N}$, where $f$ is defined by
$f(a)=a^{-1}b$, $f(b)=a^{-1}$. Let $E$ be the semidirect product
of $N$ and $K$ via the natural action of $K$ on $N$. We claim that
$E$ is non-abelian, and isomorphic to $\CoxGroup{n}$ with
generators $x=(1,\,f)$ and $y=(a,\,f)$. We have
$(a,\,f)^2=(af(a),\,f^2)=(b,\,f^2)$,
$(b,\,f^2)(1,\,f)=(b,\,\mathrm{id})$, and
$(1,\,f)(b,\,f^2)=(a^{-1},\,\mathrm{id})$. Thus $E$ is
non-abelian, and generated by $(1,\,f)$, $(a,\,f)$. A routine
computation shows that
$(1,\,f)^3=(a,\,f)^3=((1,\,f)(a,\,f))^3=((1,\,f)^{-1}(a,\,f))^n=1$.

The group $E$ proves that $|G|=3|H|=3n^2$. In particular,
$H\cong\cyclic{n}\times\cyclic{n}$.
\end{proof}
\end{theorem}

We would like to give a detailed description of the lattice of subgroups of
$\CoxGroup{p}$ in terms of generators $x$ and $y$. From a group-theoretical
point of view, the groups are rather boring, nevertheless, the lattice can be
nicely visualized. The cases $p=2$ and $p=3$ cause troubles, and \emph{we
exclude them from our discussion for the time being}.

\begin{lemma}
Let $G$ and $H$ be defined as before. Then $H$ is the Sylow
$p$-subgroup of $G$, and contains $p+1$ subgroups
$H(i)=\span{h(i)}$, for $0\leqslant i<p$, or $p=\infty$, all
isomorphic to $\cyclic{p}$. We can take
\begin{displaymath}
    h(i)=x^2y(xy^2)^i,\;\,\textrm{for \ $0\leqslant i<p$ \ and \ } h(\infty)=xy^2.
\end{displaymath}
There are $p^2$ Sylow $3$-subgroups $G(k,\,l)=\span{g(k,\,l)}$,
for $0\leqslant k$, $l<p$, all isomorphic to $\cyclic{3}$. We can
take
\begin{displaymath}
    g(k,\,l)=(x^2y)^{-k}(xy^2)^{-l}x(x^2y)^k(xy^2)^l.
\end{displaymath}
\begin{proof}
The subgroup structure of $H$ is obvious. Every element of $G\setminus H$ has
order $3$, so there are $p^2$ Sylow $3$-subgroups of order $3$ in $G$. The
subgroup $H$ acts transitively on the set of Sylow $3$-subgroups. (By Sylow
Theorems, $G$ acts transitively on the copies of $\cyclic{3}$. As $|G|=3p^2$,
the stabilizer of each $\cyclic{3}$ under this action is isomorphic to
$\cyclic{3}$. Since $p$ and $3$ are relatively prime, no element of $H$ can be
found in any stabilizer.) This shows that our list of Sylow $3$-subgroups is
without repetitions, thus complete.
\end{proof}
\end{lemma}

For certain values of $p$ (see below), there are no other subgroups in $G$.
For the remaining values of $p$, there are additional subgroups of order $3p$.

If $K\leqslant G$ has order $3p$, it contains a unique normal
subgroup of order $p$, say $L\leqslant H$. Since $L$ is normalized
by both $K$ and $H$, it is normal in $G$. Then $G/L$ is a
non-abelian group of order $3p$, and has therefore $p$ subgroups
of order $3$. Using the correspondence of lattices, we find $p$
subgroups of order $3p$ containing $L$ (the group $K$ is one of
them).

\begin{lemma}
The group $H(i)$ is normal in $G$ if and only if
\begin{equation}\label{Eq:QCongruence}
    i^2+i+1\equiv 0\,({\rm mod}\,p).
\end{equation}
If $p\equiv 1\pmod{3}$, there are two solutions to $(\ref{Eq:QCongruence})$.
For other values of $p$, there is no solution.
\begin{proof}
We have
\[\arraycolsep=.5mm\begin{array}{rcl}
    x^{-1}h(i)x&=&x^{-1}x^2y(xy^2)^ix=xy^2y^2(xy^2)^ix\\[1mm]
        &=&(xy^2)(y^2x)^{i+1}=(x^2y)^{-(i+1)}(xy^2).
\end{array}
\]
Thus $x^{-1}h(i)x$ belongs to $H(i)$ if and only if
$(x^2y)^{-(i+1)i}(xy^2)^i=(x^2y)(xy^2)^i$, i.e.  if and only if
$i$ satisfies $(\ref{Eq:QCongruence})$. Similarly,
\[\arraycolsep=.5mm\begin{array}{rcl}
    y^{-1}h(i)y&=&y^{-1}x^2y(xy^2)^iy=(y^2x)(xy^2)y^2(xy^2)^iy\\[1mm]
        &=&(y^2x)(xy^2)(y^2x)^i=(x^2y)^{-(i+1)}(xy^2).
\end{array}
\]
Then $y^{-1}h(i)y$ belongs to $H(i)$ if and only if $i$ satisfies
$(\ref{Eq:QCongruence})$.

The quadratic congruence $(\ref{Eq:QCongruence})$ has either two solutions or
none. Pick $a\in\field{p}^*$, $a\ne 1$. Then $a^2+a+1=0$ if and only if
$a^3=1$, since $a^3-1=(a-1)(a^2+a+1)$. This simple argument shows that
$(\ref{Eq:QCongruence})$ has a solution if and only if $3$ divides
$p-1=|\field{p}^*|$.
\end{proof}
\end{lemma}

\begin{theorem}[The Lattice of Subgroups of $\CoxGroup{p}$]
\label{Th:SubgroupStructure} For a prime $p>3$, let
$G=\CoxGroup{p}$, $H=\span{x^2y,\,xy^2}$, $h(i)=x^2y(xy^2)^i$ for
$0\leqslant i<p$, $h(\infty)=xy^2$, $H(i)=\span{h(i)}$,
$g(k,\,l)=(x^2y)^{-k}(xy^2)^{-l}x(x^2y)^k(xy^2)^l$ for $0\leqslant
k$, $l<p$, and $G(k,\,l)=\span{g(k,\,l)}$.

Then $H(\infty)\cong\cyclic{p}$, $H(i)\cong\cyclic{p}$,
$G(k,\,l)\cong\cyclic{3}$ are the minimal subgroups of $G$, and
$H(i)\join H(j)=H\cong\cyclic{p}\times\cyclic{p}$ for every $i\ne
j$. When $3$ does not divide $p-1$, there are no other subgroups
in $G$. Otherwise, there are additional $2p$ non-abelian maximal
subgroups of order $3p$; $p$ for each $1<i<p$ satisfying
$i^3\equiv 1\pmod{p}$. These subgroups can be listed as
$K(i,\,l)=H(i)\join G(0,\,l)$, for $0\leqslant l<p$. Then
$H(i)\join G(k',\,l')=K(i,\,l)$ if and only if $l'-l\equiv
ik'\pmod{p}$; otherwise $H(i)\join G(k',\,l')=G$. Finally,  let
$(k,\,l)\ne(k',\,l')$. Then $G(k,\,l)\join G(k',\,l')=H(i)\join
G(k,\,l)$ if and only if there is $1<i<p$ satisfying $i^3\equiv
1\pmod{p}$ such that $l'-l\equiv (k'-k)i\pmod{p}$; otherwise
$G(k,\,l)\join G(k',\,l')=G$.

The group $\CoxGroup{2}$ is isomorphic to $A_4$, the alternating group on $4$
points, and $\CoxGroup{3}$ is the unique non-abelian group of order $27$ and
exponent $3$.
\begin{proof}
Check that $h(i)^{-1}g(k,\,l)h(i)=g(k+1,\,l+i)$, and conclude that $H(i)\join
G(k,\,l)=H(i)\join G(k',\,l')$ if and only if $l'-l\equiv i(k'-k)\pmod{p}$.
This also implies that, for some $1<i<p$, $H(i)\join G(k',\,l')$ equals
$K(i,\,l)$ if and only if $l'-l\equiv ik'\pmod{p}$ and $i^3\equiv 1\pmod{p}$.

Finally, if $S=G(k,\,l)\join G(k',\,l')\ne G$, it contains a unique
$H(i)\unlhd G$. Moreover, we have $S=H(i)\join G(k,\,l)=H(i)\join G(k',\,l')$
solely on the grounds of cardinality, and everything follows.
\end{proof}
\end{theorem}

We illustrate Theorem \ref{Th:SubgroupStructure} with $p=7$. The
congruence $(\ref{Eq:QCongruence})$ has two solutions, $i=2$ and
$i=4$. The subgroup lattice of $\CoxGroup{7}$ is depicted in the
$3$D Figure \ref{Fg:SubgroupLattice}. The $49$ subgroups
$G(k,\,l)$ are represented by a parallelogram that is thought to
be in a horizontal position. All lines connecting the subgroups
$G(k,\,l)$ with $K(2,\,0)$ and $K(4,\,0)$ are drawn. The lines
connecting the subgroups $G(k,\,l)$ with $K(2,\,j)$, $K(4,\,j)$,
for $1\leqslant j< p$, are omitted for the sake of transparency.
The best way to add these missing lines is by the means of affine
geometry of $\field{p}\times\field{p}$. To determine which groups
$G(k,\,l)$ are connected to the group $K(i,\,j)$, start at
$G(0,\,j)$ and follow the line with slope $i$, drawn modulo the
parallelogram.

\setlength{\unitlength}{1.0mm}
\begin{figure}
    \centering
    \input{Coxeter.lp}
    \caption[]{The lattice of subgroups of $\CoxGroup{7}$}
    \label{Fg:SubgroupLattice}
\end{figure}

The group $A_4$ fits the description of Theorem \ref{Th:SubgroupStructure},
too, as can be seen from its lattice of subgroups in Figure \ref{Fg:A4}. So
does the group $\CoxGroup{3}$.

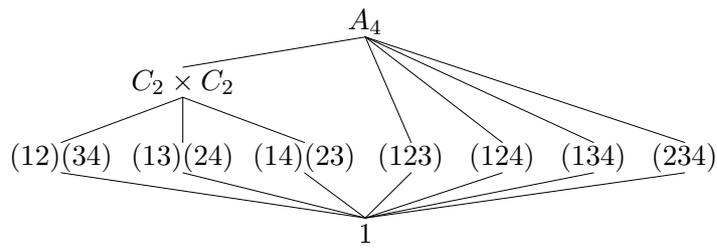
\placedrawing{a4.lp}{The subgroup structure of $A_4$}{Fg:A4}

\section*{\centerline{3. Three subgroups}}

\noindent We promised to show that each of the subgroups
$\span{g_3,\,g_4}$, $\span{g_3,\,g_5}$, $\span{g_4,\,g_5}$ of
$\paige{q}$ is isomorphic to $\CoxGroup{p}$.
\\

\noindent
{\bf Proposition 3.1.} {\em Let $g_3$, $g_4$, $g_5$ be
defined as above, $q=p^r$. Then the three subgroups $\span{g_3,\,
g_4}$, $\span{g_3,\, g_5}$, $\span{g_4,\, g_5}$ of $\paige{p^r}$
are isomorphic to $\CoxGroup{p}$, if $q\ne 9$ is odd or $q=2$.}

\begin{proof}
We prove that $G_1=\span{g_3,\,g_4}\cong\CoxGroup{p}$; the
argument for the other two groups is similar. We have
$g_3^3=g_4^3=(g_3g_4)^3=(g_4g_3)^3=(g_3^{-1}g_4)^p=(g_3^2g_4)^p=e$.
Thus $G_1\leqslant\CoxGroup{p}$. Also, $H_1=\span{g_3^2g_4,\,
g_3g_4^2}\cong\cyclic{p}\times\cyclic{p}$. When $p\ne 3$, we
conclude that $|G_1|=3p^2$, since $G_1$ contains an element of
order $3$. When $p=3$, we check that $g_3\not\in H_1$, and reach
the same conclusion.
\end{proof}

We finish this paper with a now obvious observation, that in order to describe
all subloops of $\paige{q}$, one only has to study the interplay of the
isomorphic subgroups $\span{g_3,\,g_4}$, $\span{g_3,\,g_5}$, and
$\span{g_4,\,g_5}$.

\vspace{5mm}
 \footnotesize{Department of Mathematics\hfill Received
May 7, 2001

Iowa State University

Ames, IA 50011

U.S.A.

\vspace{1mm}
petr@iastate.edu
 }

\end{document}